# REJOINDER: ONE-STEP SPARSE ESTIMATES IN NONCONCAVE PENALIZED LIKELIHOOD MODELS

By Hui Zou[1] and Runze Li[2]

*University of Minnesota and Pennsylvania State University*

We would like to take this opportunity to thank the discussants for their thoughtful comments and encouragements on our work. The discussants raised a number of issues from theoretical as well as computational perspectives. Our rejoinder will try to provide some insights into these issues and address specific questions asked by the discussants.

Most traditional variable selection criteria, such as the AIC and the BIC, are (or are asymptotically equivalent to) the penalized likelihood with the $L_0$ penalty, namely, $p_\lambda(|\beta|) = \frac{1}{2}\lambda^2 I(|\beta| \neq 0)$, and with appropriate values of $\lambda$ (Fan and Li [7]). In general, the optimization of the $L_0$-penalized likelihood function via exhaustive search over all subset models is an NP-hard computational problem. Donoho and Huo [3] and Donoho and Elad [2] show that, under some conditions, the solution to the $L_0$-penalized problem can be found by solving a convex optimization problem of minimizing the $L_1$-norm of the coefficients, when the solution is sufficiently sparse. In other words, the NP-hard best subset variable selection can be solved by efficient convex optimization algorithms under the sparsity assumption. This sheds light on variable selection for high-dimensional models, and motivates us to use continuous penalties, such as the $L_1$ penalty, rather than the discontinuous penalties, including the $L_0$ penalty. The penalized likelihood procedure with the $L_1$ penalty coincides with LASSO. In the same spirit of LASSO, the penalized likelihood with a nonconcave penalty, such as the SCAD penalty, has been proposed for variable selection in Fan and Li [5]. LASSO and the SCAD represent the two main streams of penalization method for variable selection in the recent literature. Although both methods operate as continuous thresholding rules, they appear to be very different theoretically and

Received November 2007; revised November 2007.
[1]Supported by NSF Grant DMS-07-06733.
[2]Supported by NSF Grant DMS-03-48869 and partially supported by National Institute on Drug Abuse (NIDA) Grant P50 DA10075.

computationally. The SCAD is asymptotically unbiased and enjoys the oracle properties which LASSO cannot possess (Zou [16]). On the other hand, LASSO uses the $L_1$ penalty as opposed to the concave SCAD penalty, thus being computationally more friendly than the SCAD.

Computational efficiency of a statistical method determines its popularity to a great extent. Interestingly, this argument is evidenced by the evolution of LASSO. Although the original LASSO was shown to be a promising variable selection method in Tibshirani [13], it did not become popular in statistical practice before 2002 due to the relative inefficiency of the original LASSO algorithm (Madigan and Greg [11]). The situation dramatically changed in 2002 when Efron et al. [4] invented the LARS algorithm which can compute the entire LASSO solution paths in a very efficient fashion. Since then, LASSO has enjoyed its enormous popularity.

The LQA algorithm by Fan and Li [5] provides a unified view of LASSO and the SCAD by treating them as iteratively re-weighted ridge regression. However, the unification based on LQA is not very satisfactory for two reasons. First, LASSO enthusiasts may not accept the LQA algorithm as the *algorithm* for computing LASSO, given the extremely successful LARS algorithm. Second and more importantly, re-weighted ridge regression cannot automatically produce a sparse solution. One has to fill the conceptual gap by artificially thresholding small values to zero within each iteration. As discussed in our paper, this practice is not very satisfactory.

The LLA algorithm presented in our work provides a better unification of LASSO and the SCAD, because within each iteration the solution naturally adopts a sparse representation. Putting LASSO in the LLA framework, we see that the solution converges after the first step. For the SCAD, although multiple steps are needed before convergence, the one-step estimates work as well as the final estimates as long as the initial solution is a root-$n$ consistent estimator. We have shown that the LARS algorithm can be exploited to solve the one-step sparse estimator. Considering that the one-step sparse estimator enjoys the oracle properties, we hope our results would persuade more people to use the one-step sparse estimation idea.

If one accepts the unification of LASSO and the SCAD (and many other nonconcave penalized methods), then it is helpful to bring the adaptive LASSO (Zou [16]) into the big picture. The adaptive LASSO was proposed to fix several theoretical drawbacks of LASSO. Meinshausen and Bühlmann [12] showed that LASSO selection can be consistent under the neighborhood stability condition. Later, Zou [16] derived a necessary condition for LASSO selection to be consistent and the necessary condition was shown to be equivalent to the neighborhood stability condition (Zhao and Yu [15]). However, the neighborhood stability condition is quite strong and difficult to check in practice. We could only hope the condition holds when applying LASSO to do variable selection. Furthermore, LASSO cannot have the oracle properties of the SCAD. To overcome these difficulties, Zou [16] proposed using the adaptively weighted $L_1$ penalty to replace the $L_1$ penalty in LASSO, and the modified LASSO was named the adaptive LASSO. Consider the penalized least squares. The adaptive LASSO is formulated as follows:

$$\widehat{\beta}(\text{AdaLasso}) = \arg\min_{\beta} \|\mathbf{y} - \mathbf{X}\beta\|^2 + \sum_{j=1}^{p} |\hat{\beta}_j^0|^{-\gamma} |\beta_j|,$$

where $\widehat{\beta}^0$ is a root-$n$ estimator. Zou [16] suggested to pick $\gamma$ from $\{0.5, 1, 2\}$ by cross-validation. As pointed out by Bühlmann and Meier, we could view $\widehat{\beta}(\text{AdaLasso})$ as the one-step estimator for some Type 1 penalty function. For instance, for $\gamma = 0.5$, the penalty function is the Bridge penalty $L_{0.5}$ and for $\gamma = 1$, the penalty function is the log-penalty. Thus, we happily see that the LLA framework tightly connects the three sparse estimation techniques. Based on our experience, often the adaptive LASSO with $\gamma \geq 1$ gives the best performance. However, we may feel uncomfortable to regard the adaptive LASSO as the one-step estimator when $\gamma > 1$, because the corresponding Type 1 penalty should be $p_\lambda(t) = \frac{\lambda}{1-\gamma}|t|^{1-\gamma}$, which is negative! We have not seen any use of a negative penalty function in the literature. It is hard to imagine that one would try to solve a negatively penalized least-squares problem in the first place.

Both Meng and Bühlmann and Meier have suggested the possibility to go beyond the one-step estimator. Meng worries about the statement "provided that the initial estimators are reasonably good" and suggests that multiple-step estimators provide some safety-net for guarding against accidental "unreasonable" starting points. We in general agree with Meng. However, we would also like to point out that, unlike in some general-purpose iterative optimization algorithms, we do often have a good initial estimator in the context of one-step sparse estimation. Although in the paper we only considered using the ordinary MLE as the initial estimator for simplicity, practically, some regularized estimators such as LASSO or the Elastic Net could yield better one-step estimators. In fact, Bühlmann and Meier have demonstrated the very competitive performance of the one-step estimator with LASSO as the initial estimates. Meng provides some Bayesian-flavored explanation to the improvement of the one-step estimator over the final estimator by making a connection to a similar phenomenon appearing in the EM literature. The analogy is quite natural since the LLA algorithm is shown to be equivalent to an EM algorithm in our paper. Here we try to understand the outperformance phenomenon from a different angle. First of all, it is easy to see that the one-step estimator would not outperform the fully-iterated estimator if the initial estimator was not good. Second, once a good initial estimator is used in the LLA algorithm, our theory says that the one-step estimator achieves the asymptotic efficiency of the final estimator. After that, the iteration could add extra variability into the estimator and eventually increase the overall mean squared error. This argument is not very precise, but similar conclusions have been made for other iterative optimization algorithms such as boosting. Boosting minimizes an empirical loss function via iterative functional gradient decent. It is now well known that boosting forever can increase the misclassification error (Hastie, Tibshirani and Friedman [9]), and thus, early stopping is necessary to avoid over-fitting (Zhang and Yu [14]). Using the one-step estimator in the LLA algorithm is similar to using the early-stopping rule in boosting.

Bühlmann and Meier suggest a new multiple-step estimator named MSA-LASSO. Their proposal is very interesting as demonstrated in their simulation study. We think MSA-LASSO can be regarded as an adaptive way to search for the best initial estimator in the one-step paradigm. What is interesting is that the one-step estimation idea itself is used in the adaptive search. We agree with them that reducing the number of false positives is perhaps more important than reducing prediction errors in high-dimensional data analysis. Thus, the simulation study indicates the promising potential of MSA-LASSO in real applications.

Fan and Li [5] and Fan and Peng [8] established the existence of a local minimizer of the nonconcave penalized likelihood which holds the optimal oracle properties. However, what is left to be shown is whether the final estimator from the LQA algorithm always finds the desired local minimizer. This issue does not exist any more in the one-step sparse estimation and the LLA algorithm. Zhang considers another innovative way to bypass the local minimizer issue. Zhang's idea is to modify the SCAD function such that the new penalty function retains the shape of the SCAD penalty but has the smallest maximum concavity. By controlling the maximum concavity, one could establish some convexity conditions to guarantee the convexity of the optimization criterion in the region of interests. For example, if we are only interested in a sparse model with at most $d^*/2$ variables, then Zhang shows that the sparse convexity condition (3.2) guarantees the solution is the unique local minimizer. We think Zhang's MC+ method provides deep and new insights into the nonconcave penalized models. In addition, Zhang also proposes the PLUS algorithm to efficiently compute the unique local minimizer. If viewing the MCP as a Type 2 penalty, we can use the LLA algorithm to compute the MC+ solution for a fixed penalization parameter. Zhang shows that the PLUS algorithm can compute the solution paths, which is a remarkable achievement in our opinion. Overall, both Zhang's idea and ours share an important common theme, that is, effectively exploiting convexity in maximizing nonconcave penalized likelihood. We agree with Bühlman and Meier that Type 1 penalty functions may enjoy better computational efficiency in the one-step estimation paradigm, however, as demonstrated by Zhang, Type 2 penalty functions may enjoy some nice properties.

In our paper we have presented a hierarchical Bayesian construction in which the corresponding EM algorithm is exactly the LLA algorithm for a large class of concave penalty functions. Our original goal was to use this result to establish the statistical foundation for the LLA algorithm. However, as outsiders (to the EM community), we did not further explore the LLA = EM result. Further research on this topic is certainly of interest. We agree with Meng that it is somewhat careless to jump to the conclusion based on Theorem 3 that MM algorithms are more flexible than EM algorithms. It is very possible that more skilled statisticians can recast the LLA algorithm for the SCAD as an EM algorithm by a cleverer construction.

On the other hand, we are glad to see our ignorance triggered some very intriguing discussion on a Bayesian approach to variable selection that is connected to but different from the penalized likelihood approach. The connection is nicely summarized in Meng's sentence: "using penalized likelihood enjoys the Bayesian fruits without paying the B-club fee." In the B-club, ridge regression and LASSO are the MAPS with Gaussian and Laplace priors, respectively. However, pena-likelihoodists also bring some new fruit to the B-club. For example, our Bayesian friends may regard the SCAD as the MAP using a SCAD prior (the density is the exponential of negative SCAD penalty function). But the SCAD prior was not on the B-club menu when Fan and Li [5] invented the SCAD estimator. Meng suggests that a classical Bayesian strategy to attack the variable selection problem would start with a prior which puts nontrivial mass at zero. The mixture formulation of such a prior has been studied in various papers. Meng also directs our attention to using the *posterior median* instead of the usual posterior mode to achieve sparsity. We find his arguments are quite intriguing and eye-opening. Barbieri and Berger [1] have shown that under the Bayesian approach, the optimal prediction model is the *median probability model*, and is not always the model with highest posterior probability. Two papers deliver strikingly similar messages to advocate the use of posterior median, although the meaning of posterior median is not exactly the same in two papers. Nevertheless, these results seem to indicate that the penalized likelihood solution may not necessarily correspond to the optimal Bayesian solution.

Finally, we agree with Meng that the Bayesian approach allows us to gain probabilistic insights via the full probabilistic modeling. Without paying the B-club fee and thus being outsiders, we seem to have difficulty in specifying a meaningful prior for some classes of semiparametric regression models. The penalized likelihood approach can be naturally extended for the partially linear models and generalized varying coefficient partially linear models, two classes of popular semiparametric regression models (Fan and Li [6] and Li and Liang [10]), although here we are not arguing that the penalized likelihood is intrinsically more flexible than the Bayesian approach. Furthermore, in many applications, such as the real data example presented by Bühlmann and Meier, it is very likely that the penalized likelihood is computationally more appealing than the existing Bayesian variable selection procedures.

**Acknowledgments.** We thank all the discussants for their stimulating comments. Our thanks also go to the Co-Editors for organizing this discussion so that we can have this nice platform to exchange ideas.

| | |
|---|---|
| School of Statistics<br>University of Minnesota<br>Minneapolis, Minnesota 55455<br>USA<br>E-mail: hzou@stat.umn.edu | Department of Statistics<br>Pennsylvania State University<br>University Park, Pennsylvania 16802–2111<br>USA<br>E-mail: rli@stat.psu.edu |